\documentclass[11pt]{amsart}
\usepackage{amsmath,amstext,amssymb,amsopn,amsthm}
\usepackage[utf8]{inputenc}
\usepackage[T1]{fontenc}

\usepackage{amscd}
\usepackage{amsaddr}
\usepackage{amsfonts}

\usepackage{graphicx}
\usepackage{color}
\usepackage{enumerate}
\usepackage{tikz}
\usepackage{float}
\usepackage[font=scriptsize]{caption}
\usepackage{verbatim}
\usepackage{dsfont}
\usepackage{hyperref}
\title{Burkholder
inequality by Bregman divergence}
\author[K. Bogdan]{Krzysztof Bogdan}
\address{Faculty of Pure and Applied Mathematics, Wroc\l aw University
	of Science and Technology, Wybrze\.ze Wyspia\'nskiego 27, 50-370
	Wroc\l aw, Poland \\ (e-mail: Krzysztof.Bogdan@pwr.edu.pl)}
\author[M. Wi\k{e}cek]{Mateusz Wi\k{e}cek}
\address{Faculty of Pure and Applied Mathematics, Wroc\l aw University
	of Science and Technology, Wybrze\.ze Wyspia\'nskiego 27, 50-370
	Wroc\l aw, Poland\\ (e-mail: Mateusz.Wiecek@pwr.edu.pl)}
	\thanks{K.B. was supported by the National Science Centre grant 2017/27/B/ST1/01339. M.W. was supported by the National Science Centre grant 2018/31/B/ST1/03818. }
\newtheorem{theorem}{Theorem}

\theoremstyle{definition}

\usepackage[outer=3cm, inner=3cm, top=1in, bottom=1.25in, marginparsep=0pt, marginparwidth=12pt]{geometry}
\newcommand{\E}{\ensuremath{\mathbb{E}\, }}
\newcommand{\N}{\ensuremath{\mathbb{N}}}
\newcommand{\R}{\ensuremath{\mathbb{R}}}
\newcommand\Snp[1][p]{\ensuremath{S_{n}^{#1}}}
\newcommand\Sjp[1][j]{\ensuremath{S_{#1}^{p}}}
\newcommand\Stp[1][j]{\ensuremath{S_{#1}^{2}}}

\begin{document}
\maketitle
\date{\today}
\begin{abstract}
We prove Burkholder inequality using Bregman divergence.
\end{abstract}
\section{Introduction}
Here is the celebrated Burkholder inequality.
\begin{theorem}[Burkholder]\label{B}
	For every $p\in(1,\infty)$ there exist constants $0<c_p\le C_p<\infty$
 %\in (0,\infty)$ 
	such that for every real-valued martingale $(X_n)_{n=1}^{\infty}$,
\begin{equation}\label{wz1}
	c_p^p\E X_n^{*p}\le \E S_n(X)^p\le C_p^p\E X_n^{*p}, \quad n=1,2,...,
\end{equation}	
where $X_n^*=\max\limits_{1\le j\le n}|X_j|$, $S_n(X)=
\bigl(\sum_{j=1}^{n}(\Delta X_j)^2\bigr)^{1/2}$, 
$\Delta X_n=X_n-X_{n-1}$ and $X_0=0$.
\end{theorem}
The inequality~\eqref{wz1} was proved by 
Burkholder in \cite[Theorem~9]{Burkholder} by using weak $(1,1)$ estimates and Marcinkiewicz interpolation theorem. For martingales satisfying additional conditions, e.g., for continuous martingales, 
Burkholder and 
Gundy \cite{BG} extended \eqref{wz1} to the range $0<p<\infty$. 
Davis proved in \cite{Davis} that the inequality~\eqref{wz1}  holds for $p=1$ and all martingales, and so \eqref{wz1} is called the Burkholder-Davis-Gundy inequality if $p\in [1,\infty)$.
Noteworthy, the case $p>1$ in \eqref{wz1} follows from the case $p=1$, see Kallenberg \cite[proof of Theorem 26.12]{Kallenberg}, see also the historical and bibliographical remarks in 
Ba\~{n}uelos \cite{MR2928339},
Marinelli and R\"{o}ckner \cite{MR3463679} and Gusak, Kukush, Kulik, Mishura and Pilipenko \cite[Chapter~7]{MR2572942}.
We refer to Beiglb\"{o}ck and Siorpaes \cite{MR3322322} and Yaroslavtsev \cite{MR4156214} for related developments, too.
The 
inequality \eqref{wz1} is one of the cornerstones of stochastic analysis, used to bound stochastic integrals, see, e.g., 
Ikeda and 
Watanabe \cite{MR1011252}.
Historically the inequality may be 
considered as far-reaching extension of  the results of 
Paley \cite{Paley} on the Walsh~series and the results of 
Zygmund and 
Marcinkiewicz \cite{zbMATH03028604} on series with independent centered increments.

We give a new proof of 
\eqref{wz1} for $1<p<\infty$ by using the so-called Bregman divergence. The proof is elementary, as can be judged from the length of this note.
It bears some similarity to the stochastic-calculus proofs in \cite{MR3463679} 
and Ba\~nuelos \cite[Theorem~1.1]{MR816309} and to the arguments of Stein \cite[Section IV.2.1]{MR0290095}
but the Bregman divergence allows us to  handle all $p>1$ directly.

\textbf{Acknowledgements.}  
The proof was drafted in 2011 after discussions  with Rodrigo Ba\~nuelos over the book of Garsia \cite{Garsia} -- the received hospitality from Purdue University is gratefully remembered. 
Unexpectedly, the draft inspired new results on operator semigroups by Ba\~nuelos, 
Bogdan and 
Luks \cite{MR3556449} and Bogdan, Jakubowski, Lenczewska and Pietruska-Pałuba \cite{MR4372148},
and harmonic functions by
Bogdan, 
Dyda and 
Luks \cite{MR3251822} and 
Bogdan, 
Grzywny, 
Pietruska-Pałuba and 
Rutkowski \cite{bogdan2020nonlinear}. Therefore we decided to
rewrite and correct the notes -- with focus on the simple. We take this opportunity to thank the authors of the papers and Adam Os\c{e}kowski for discussions that helped clarify the setting. We also thank an anonymous referee for corrections and valuable suggestions.

The proof of Theorem~\ref{B} is given in Section~\ref{s.pBG}.
Here we comment on the constants in \eqref{wz1}. 
Using the good-$\lambda$ inequalities, 
Burkholder \cite[Theorem 3.2]{Burkholder1973} proved that \eqref{wz1} holds
with $C_p=18p^{3/2}/(p-1)$ and $c_p=\bigl(18p^{5/2}/(p-1)^{3/2}\bigr)^{-1}$, see also 
Burkholder, 
Davis and 
Gundy \cite{MR0400380}.
It is known  that $c_p=1/p$ is optimal for $p\ge 2$ (see, \cite[Inequality (3.9)]{Burkholder1991}, \cite[Theorem 7.11]{Osekowski}).
 %, otherwise the optimal values of $c_p$ and $C_p$ are not known. 
In the course of the proof we obtain the following values:
	\begin{equation*}
	c_p=\begin{cases}
	\frac{1}{p\sqrt{2}}\left(\frac{p-1}{p}\right)^{p+\tfrac12} & \text{for }1<p<2, \\
	1/2 & \text{for } p=2,\\
	\frac{\sqrt 2}{p}\left(\frac{p-1}{p}\right)^{\tfrac{p-1}2} & \text{for } p>2,
	\end{cases}
	 \quad \text { and }\quad C_p=\begin{cases}
	\frac1{\sqrt{p(p-1)/2}} & \text{for } 1<p<2, \\
	1 & \text{for } p=2,\\
	\sqrt{2p} & \text{for } p>2.
	\end{cases}
	\end{equation*}
Note that the value $C_p=\sqrt{2p}$ for $p>2$ was first given by 
Garsia \cite[Theorem II.1.1]
{Garsia}. Our $C_p=\frac1{\sqrt{p(p-1)/2}}$
is better than the constant 
$(2+\sqrt{p+4})/\sqrt{p}$
from
\cite[p. 40]{Garsia} for $p\in (p_0,2)$, where $p_0\sim 1.11018$, but is worse for $p\in (1,p_0)$. 
Our value of $c_p$ 
for $1<p<2$ is worse than 
$1/\sqrt{10p}$ given in \cite[II.2.8, p. 37]{Garsia}
and  is worse for $p>2$ than the optimal value $1/p$ by Burkholder. 
On the other hand our values 
of $c_p$ and $C_p$ improve
those resulting
from the good-$\lambda$ inequalities, given above.

For $k\in \R$
we let $$x^{\langle k \rangle}=|x|^k\mathrm{sign}(x), \quad x\in\R,$$
with the conventions $0^{\langle k \rangle}=\mathrm{sign}(0)=0$. Clearly, $(x^{\langle k \rangle})'=k|x|^{k-1}$ and $(|x|^k)'=kx^{\langle k-1 \rangle}$ for $x\neq 0$. 
We define the Bregman divergence:
\begin{equation}\label{F_p}
F_p(a,b):=|b|^p-|a|^p-p(b-a)a^{\langle p-1 \rangle}, \ a,b\in\R.
\end{equation}
For instance, $F_2(a,b)=(b-a)^2$. For applications of Bregman divergence in analysis, statistics learning and optimization, we refer the reader to \cite{bogdan2020nonlinear}.
Observe that by the convexity of the function $\R\ni x\mapsto |x|^p$, we have $F_p(a,b)\ge 0$. In fact, by
\cite[Lemma 6]{MR3251822},
	for every $p\in(1,\infty)$ there are constants $d_p,D_p>0$ such that
	\begin{equation}\label{FvsG}
	d_p(b-a)^2(|a|\vee |b|)^{p-2}\le F_p(a,b)\le D_p(b-a)^2(|a|\vee |b|)^{p-2}, \quad a,b\in\mathbb{R},
	\end{equation}
with the convention that $(b-a)^2(|a|\vee |b|)^{p-2}=0$ if $a=b=0$.	
We claim that the values $D_p=p(p-1)/2$ for $p\ge2$ and $d_p=p(p-1)/2$ for $p\in(1,2)$ suffice and are optimal in \eqref{FvsG}. 
Indeed, we first check the trivial case $a=0$. Then we assume $a\neq 0$ and denote $G_p(a,b):=(b-a)^2(|a|\vee |b|)^{p-2}$. We get
$F_p(a,b)=|a|^pF_p(1,b/a)$ and $G_p(a,b)=|a|^pG_p(1,b/a)$, so it is enough to consider $a=1$.
We have
\begin{equation}\label{integrals}
F_p(1,b)=p(p-1)\int_{1}^b\int_1^x|y|^{p-2}\mathrm{d}y\mathrm{d}x
\end{equation}
and, of course,
\begin{equation*}\label{integralsG}
G_p(1,b)=2\int_{1}^b\int_1^x(|b|\vee 1)^{p-2}\mathrm{d}y\mathrm{d}x.
\end{equation*}
The claim follows since  in \eqref{integrals} we have $|y|\le |b|\vee 1$ and $|y|\sim |b|\vee 1\sim 1$ when $b\sim 1$. Our explicit values of $c_p$ and $C_p$ above depend on this result since in the proof of  \eqref{wz1} below we will use $d_p$ only
when $1<p< 2$ and 
$D_p$ only when $p> 2$.

 Let $f,g$ be two nonnegative functions with the same domain. We say that $f$ and $g$ are comparable and write $f\asymp g$, if there exist constants $c,C>0$ such that 
	\begin{equation}
	cf(x)\le g(x)\le Cf(x)
	\end{equation}
for all arguments~$x$.	
Thus, \eqref{FvsG} reads as $F_p(a,b)\asymp (b-a)^2(|a|\vee |b|)^{p-2}$. 

\section{Proof of the Burkholder inequality}\label{s.pBG}	
 Let $X=(X_n)_{n=1}^{\infty}$ be adapted to the filtration $(\Omega, (\mathcal{F}_n)_{n=1}^{\infty}, \mathbb{P})$. 
In addition to the notation in Theorem~\ref{B}, we denote $S_0(X)=0$
and
$S_n^p(X)=(S_n(X))^p$.
Doob's inequality \cite[Theorem~3.4, Chapter VII]{Doob} yields
	\begin{equation}\label{Doob_inequalities}
	\E|X_n|^p\le \E(X_n^*)^p\le \Bigl(\frac{p}{p-1}\Bigr)^p\E|X_n|^p, \ n=1,2,...\,, \quad p\in(1,\infty) .
	\end{equation}

Let
$n\in\N$ and 
$p\in(1,\infty)$. Of course, if $\E (X_n^*)^p=0$, then $X_j=0$ a.s. for $j=1,...,n$, so $\E \Snp(X)=0$ and the inequalities \eqref{wz1} obviously hold with any positive values of $c_p,C_p$.  Similarily, if $\E|X_j|^p=\infty$ and $j\le n$, then $\E(X_n^*)^p=\infty$ and $\E S_n^p(X)=\infty$, so we have \eqref{wz1} with any $c_p,C_p>0$. Therefore in what follows we assume that $\E (X_n^*)^p>0$ and 
$\E |X_j|^p<\infty$ for $j=1,2,...,n$, i.e., $\E(X_n^*)^p<\infty$. 
	Clearly,
	\begin{equation}\label{wz11}
	|X_n|^p=\sum_{j=1}^{n}(|X_j|^p-|X_{j-1}|^p).
	\end{equation}
	Observe that if $X_n\in L^p:=L^p(\mathbb{P})$, then $X_n^{\langle p-1 \rangle}\in L^{q}(\mathbb P)$ for $q=p/(p-1)$. In fact, $\E|X_n^{\langle p-1 \rangle}|^{p/(p-1)}=\E|X_n|^p$. Of course, $\Delta X_j\in L^p$ for every $j=1,2,...,n$. We get
	\begin{align*}
	p\E[(X_j-X_{j-1})X_{j-1}^{\langle p-1 \rangle}]&=p\E(\E[\Delta X_jX_{j-1}^{\langle p-1 \rangle}|\mathcal{F}_{j-1}])\\&=p\E(X_{j-1}^{\langle p-1 \rangle}\E[\Delta X_j|\mathcal{F}_{j-1}])=0.
	\end{align*}
	The expectations above are absolutely convergent by H{\"o}lder's inequality (such 
explanations are usually omitted later on).
	From \eqref{wz11} it follows that
	\begin{equation}\label{F_przydatne}
	\E|X_n|^p=\E\sum_{j=1}^{n}F_p(X_{j-1},X_j).
	\end{equation}

The case $\mathbf{p=2}$ in \eqref{wz1} is rather trivial.
Indeed, $F_2(a,b)=(b-a)^2$. By \eqref{F_przydatne}, 
$\E |X_n|^2= \E S_n^2(X)$, so by Doob's inequality, 
	\begin{equation}
	\frac{1}{4}\E(X_n^*)^2\le \E S_n^2(X)\le \E (X_n^*)^2.
	\end{equation}
We next consider the case $\mathbf{p>2}$.
By \eqref{FvsG}  and~\eqref{F_przydatne},
	\begin{equation}\label{F_przydatne2}
	\E|X_n|^p\le D_p\E\Bigl(\sum_{j=1}^n(X_{j}-X_{j-1})^2(|X_j|\vee |X_{j-1}|)^{p-2}\Bigr).
	\end{equation}
	By 
\eqref{F_przydatne2} and Doob's inequality, 
	\begin{equation}\label{Cos_z_Doobem}
	\left(\frac{p-1}{p}\right)^p\E (X_n^*)^p\le \E|X_n|^p\le D_p\E(X_n^*)^{p-2}\Snp[2](X).
	\end{equation}
	Using H{\"o}lder's inequality with exponents $p/(p-2)$
and $p/2$
we get
	\begin{equation*}
	\E(X_n^*)^{p-2}\Snp[2](X)\le\left(\E(X_n^*)^p\right)^{\frac{p-2}{p}}\left(\E \Snp(X)\right)^{\frac{2}{p}}.
	\end{equation*} 
By  \eqref{Cos_z_Doobem},
	\begin{equation*}
	\left(\frac{p-1}{p}\right)^p\E (X_n^*)^p\le D_p\left(\E(X_n^*)^p\right)^{\frac{p-2}{p}}\left(\E \Snp(X)\right)^{\frac{2}{p}}.
	\end{equation*}
	Dividing by $\left(\E(X_n^*)^p\right)^{\frac{p-2}{p}}$ (this is allowed since $\E(X_n^*)^p>0$), we obtain
	\begin{equation*}
		\left(\frac{p-1}{p}\right)^p\left(\E (X_n^*)^p\right)^{\frac{2}{p}}\le D_p\left(\E \Snp(X)\right)^{\frac{2}{p}},
	\end{equation*}
	or, equivalently,
	\begin{equation}\label{p2_1}
	{D_p}^{-p/2}\left(\frac{p}{p-1}\right)^{-p^2/2}\E (X_n^*)^p\le \E \Snp(X).
	\end{equation}
	 We next prove the reverse of \eqref{p2_1} (the argument is similar to that in 
\cite[Theorem~II.1.1,~p.~28]{Garsia}). Like in \eqref{wz11}, we have
	\begin{equation}
	\Snp(X)=\sum_{j=1}^{n}(\Sjp(X)-\Sjp[j-1](X)).
	\end{equation}
	Let $j_0=n\land \min\{1\le j\le n: S_{j}(X)>0\}$. Then,
	\begin{equation}\label{wzorek2}
	\Snp(X)=\sum_{j_0<j\le n}(\Sjp(X)-\Sjp[j-1](X))+\Sjp[j_0](X).
	\end{equation}
For $\alpha, b\ge1$
we have
	\begin{equation}\label{oczywista}
	b^{\alpha}-1=\alpha\int_{1}^{b}t^{\alpha-1}\mathrm{d}t\le \alpha b^{\alpha-1}(b-1).
	\end{equation} 
	By \eqref{wzorek2} and \eqref{oczywista} we obtain
	\begin{flalign*}
	\Snp(X)&=\sum_{j_0<j\le n}\Sjp[j-1](X)\Bigl(\Bigl(\frac{\Stp(X)}{\Stp[j-1](X)}\Bigr)^{\frac{p}{2}}-1\Bigr) +S_{j_0}^p(X) \\ &\le \sum_{j_0<j\le n}\Sjp[j-1](X)\frac{p}{2}\Bigl(\frac{\Stp(X)}{\Stp[j-1](X)}\Bigr)^{\frac{p}{2}-1}\Bigl(\frac{\Stp(X)}{\Stp[j-1](X)}-1\Bigr)+\Sjp[j_0](X)\\&\le\frac{p}{2}\Bigl(\sum_{j_0<j\le n}(\Stp(X)-\Stp[j-1](X))S_{j}^{p-2}(X)+\Sjp[j_0](X)\Bigr)\\&=\frac{p}{2}\sum_{j=1}^n(\Delta X_j)^2\sum_{l=1}^{j}\Delta S_{l}^{p-2}(X).
	\end{flalign*}
	Changing the order of summation we get
	\begin{equation*}
	\Snp(X)\le\frac{p}{2}\sum_{l=1}^{n}\Delta S_{l}^{p-2}(X)\sum_{j=l}^n(\Delta X_j)^2,
	\end{equation*}
	thus
	\begin{flalign*}
	\E\Snp(X)&\le \frac{p}{2}\E\Bigl(\sum_{l=1}^{n}\Delta S_{l}^{p-2}(X)\sum_{j=l}^{n}(X_j-X_{j-1})^2\Bigr)\\&=\frac{p}{2}\sum_{l=1}^{n}\E\Bigl(\Delta S_l^{p-2}(X)(X_l-X_{l-1})^2+\sum_{l+1\le j\le n}\E[\Delta S_{l}^{p-2}(X)(X_j-X_{j-1})^2|\mathcal{F}_{j-1}]\Bigr)\\&=\frac{p}{2}\sum_{l=1}^{n}\E\Bigl(\Delta  S_l^{p-2}(X)\Bigl((X_l-X_{l-1})^2+\sum_{l+1\le j\le n}\E[(X_j-X_{j-1})^2|\mathcal{F}_{j-1}]\Bigr)\Bigr).
	\end{flalign*}
	We have
	 \begin{flalign*}
	 \E[(X_j-X_{j-1})^2|\mathcal{F}_{j-1}]&=\E[ X_j^2|\mathcal{F}_{j-1}]+\E[ X_{j-1}^2|\mathcal{F}_{j-1}]-2\E[X_jX_{j-1}|\mathcal{F}_{j-1}]\\&=\E[ X_j^2|\mathcal{F}_{j-1}]-X_{j-1}^2=\E[ X_j^2-X_{j-1}^2|\mathcal{F}_{j-1}],
	 \end{flalign*}
	therefore,
	 \begin{flalign}
	\E\Snp(X)&\le\frac{p}{2}\sum_{l=1}^{n}\E\Bigl(\Delta S_l^{p-2}(X)\Bigl((X_l-X_{l-1})^2+\sum_{l+1\le j\le n}(X_j^2-X_{j-1}^2)\Bigr)\Bigr)\notag\\&=\frac{p}{2}\sum_{l=1}^{n}\E\Bigl(\Delta S_l^{p-2}(X)\Bigl(X_n^2+X_{l-1}^2-2X_lX_{l-1}\Bigr)\Bigr)\label{125}\\&\le 2p\E\Bigl(\sum_{l=1}^{n}\Delta{S_{l}}^{p-2}(X)(X_n^*)^2\Bigr)=2p\E[S_{n}^{p-2}(X)(X_n^*)^2].\notag
	 \end{flalign}
	 By H{\"o}lder's inequality with exponents $p/(p-2)$ and $p/2$,
	 \begin{equation}\label{126}
	 2p\E[S_{n}^{p-2}(X)(X_n^*)^2]\le 2p\left(\E S_{n}^{p}(X)\right)^{\frac{p-2}{p}}\left(\E{(X_n^*)}^{p}\right)^{\frac{2}{p}}.
	 \end{equation}
	 From \eqref{125} and \eqref{126} we obtain
	 \begin{equation}\label{p2_2}
	 \E\Snp(X)\le (2p)^{\tfrac{p}2}\E{(X_n^*)}^{p}.
	 \end{equation}
	  Combining  \eqref{p2_1} and \eqref{p2_2} with  the observation that $D_p=p(p-1)/2$ for $p>2$,
	  we get
	 \begin{equation}\label{ineq2}
	 \frac{2^{\tfrac p2}}{p^p}\left(\frac{p-1}{p}\right)^{\tfrac{p(p-1)}2}\E (X_n^*)^p\le\E\Snp(X)\le (2p)^{\tfrac{p}2}\E{(X_n^*)}^{p},
	 \end{equation}
	 or
	 \begin{equation}\label{ineq2p}
	 \frac{\sqrt 2}{p}\left(\frac{p-1}{p}\right)^{\tfrac{p-1}2}\bigl(\E (X_n^*)^p\bigr)^{\tfrac{1}p}\le\bigl(\E\Snp(X)\bigr)^{1/p}\le \sqrt{2p}\bigl(\E{(X_n^*)}^{p}\bigr)^{\tfrac{1}p}.
	 \end{equation}
It remains to consider the case $\mathbf{1<p<2}$.
Let $q=p/(p-1)$. Clearly, $q>2$. 
Assume that $Z=(Z_j)_{j=1}^{n}$ is a martingale adapted to the same filtration $(\mathcal{F}_j)_{j=1}^{n}$ and $\E |Z_j|^q<\infty$ for $j=1,...,n$. Then,
	 \begin{equation}
	 X_nZ_n=\sum_{j=1}^{n}\Delta X_j\sum_{l=1}^{n}\Delta Z_l=\sum_{j=1}^{n}\sum_{l=1}^{n}\Delta X_j \Delta Z_l.
	 \end{equation}
	 By the martingale property, $\E\Delta X_j\Delta Z_l=0$ for all $j\neq l$, hence
	 \begin{equation*}
	 \bigl|\E X_nZ_n\bigr|=\Bigl|\E\sum_{j=1}^{n}\Delta X_j \Delta Z_j\Bigr|\le \E\left(\sqrt{\sum_{j=1}^{n}(\Delta X_j)^2}\sqrt{\sum_{j=1}^{n}(\Delta Z_j)^2}\right)=\E [S_n(X)S_n(Z)].
	 \end{equation*}
	 Applying H{\"o}lder's inequality with exponents $p$ and $q$ we obtain
	 \begin{equation}\label{Holderp1}
	 |\E X_nY_n|\le \bigl(\E S_n^p(X)\bigr)^{\tfrac{1}p}\bigl(\E S_n^q(Z)\bigr)^{\tfrac{1}q}.
	 \end{equation}
	 We let $Y_n=X_n^{\langle p-1 \rangle}$. Since $|X_n|^p=|Y_n|^q$, $Y_n\in L^{q}$. Define martingale $Z=(Z_{j})_{j=1}^{n}$ by $$Z_{j}=\E[Y_n|\mathcal{F}_j],\quad  j=1,...,n.$$ Obviously, 
	 $Z_{n}=Y_n$. By Doob's inequality and \eqref{Holderp1} we get
	 \begin{equation*}
	 \Bigl(\frac{p-1}{p}\Bigr)^p\E({X_n^*})^p\le \E|X_n|^p= \E X_nZ_n\le \bigl(\E S_n^p(X)\bigr)^{\tfrac{1}p}\bigl(\E S_n^q(Z)\bigr)^{\tfrac{1}q}.
	 \end{equation*}
	 Applying \eqref{ineq2} to $\E {S_n^q}(Z)$ with $p$ replaced by $q$ we get
	 \begin{equation*}
	 \Bigl(\frac{p-1}{p}\Bigr)^p\E({X_n^*})^p\le \sqrt{2q}\bigl(\E S_n^p(X)\bigr)^{\tfrac{1}p}\bigl(\E{(Z_{n}^*)}^{q}\bigr)^{\tfrac{1}q}.
	 \end{equation*}
	 Moreover, by Doob's inequality we obtain
	 \begin{equation*}
	 \E{(Z_{n}^*)}^{q}\le\Bigl(\frac{q}{q-1}\Bigr)^q\E|Z_{n}|^q=p^{q}\E|X_n|^p\le p^{q}\E(X_n^*)^p.
	 \end{equation*}
	 Therefore,
	 \begin{equation*}
	 \Bigl(\frac{p-1}{p}\Bigr)^p\E({X_n^*})^p\le \sqrt{\frac{2p^3}{p-1}}\bigl(\E{(X_n^*)}^{p}\bigr)^{\tfrac{1}{q}}\bigl(\E S_n^p(X)\bigr)^{\tfrac{1}p},
	 \end{equation*}
	 so
	 \begin{equation*}
	 \left(\E({X_n^*})^p\right)^{\tfrac{1}p}\le p\sqrt{2}\left(\frac{p}{p-1}\right)^{p+\frac{1}{2}}\left(\E S_n^p(X)\right)^{\tfrac{1}p},
	 \end{equation*}
	 or
	 \begin{equation}\label{p1_1}
	 (2p^2)^{-p/2}\left(\frac{p-1}{p}\right)^{p\left(p+\tfrac12\right)}\E({X_n^*})^p\le \E S_n^p(X).
	 \end{equation}
	  To prove a reverse inequality, we use \eqref{F_przydatne} once more. By  \eqref{FvsG},
	 \begin{equation}\label{zdp}
	 \E|X_n|^p=\E\Bigl(\sum_{j=1}^{n}F_p(X_{j-1},X_j)\Bigr)\ge d_p\E\Bigl(\sum_{j=1}^n(X_j-X_{j-1})^2(|X_j|\vee |X_{j-1}|)^{p-2}\Bigr).
	 \end{equation}
	 Since $p-2<0$, the following inequality holds
	 \begin{flalign*}
	 \E\Bigl(\sum_{j=1}^n(X_j-X_{j-1})^2(|X_j|\vee |X_{j-1}|)^{p-2}\Bigr)&\ge \E\Bigl(\sum_{j=1}^n(X_j-X_{j-1})^2(X_n^*)^{p-2}\Bigr)\\&=\E\bigl[\Stp[n](X)(X_n^*)^{p-2}\bigr].
	 \end{flalign*}
	 By H{\"o}lder's inequality with exponents $2/(2-p)$
	 and 
$2/p$,
	 \begin{flalign*}
	 \E[S_n^p(X)]&=\E\Bigl[\Bigl(\frac{X_n^*}{S_n(X)}\Bigr)^{\frac{p(2-p)}2}\Bigl(\frac{X_n^*}{S_n(X)}\Bigr)^{\frac{p(p-2)}2}S_n^p(X)\Bigr]\\&\le \Bigl(\E\Bigl[\Bigl(\frac{X_n^*}{S_n(X)}\Bigr)^{p}S_n^p(X)\Bigr]\Bigl)^{\frac{2-p}2}\Bigl(\E\Bigl[\Bigl(\frac{X_n^*}{S_n(X)}\Bigr)^{p-2}S_n^p(X)\Bigr]\Bigl)^{\tfrac{p}2}\\&=\Bigl(\E (X_n^*)^p\Bigr)^{\frac{2-p}{2}}\Bigl(\E\bigl[\Stp[n](X)(X_n^*)^{p-2}\bigr]\Bigl)^{\tfrac{p}2}.
	 \end{flalign*}
	 Hence by \eqref{zdp} (and the assumption $\E (X_n^*)^p>0$),
	 \begin{equation*}
	 \E|X_n|^p\ge d_p \E\bigl[\Stp[n](X)(X_n^*)^{p-2}\bigr]\ge d_p\frac{\Bigl(\E\Snp(X)\Bigr)^{\tfrac{2}p}}{\Bigl(\E (X_n^*)^p\Bigr)^{\frac{2-p}{p}}}.
	 \end{equation*}
	 Since $\E (X_n^*)^p\ge \E|X_n|^p$, we get
	 \begin{equation}\label{p1_2}
	 {d_p}^{-\tfrac{p}2}\E (X_n^*)^p\ge \E\Snp(X).
	 \end{equation}

	 Combining \eqref{p1_1} and \eqref{p1_2} with the observation that $d_p=p(p-1)/2$ for $p\in (1,2)$, we get
	 \begin{equation}
	 \frac{1}{(2p^2)^{\tfrac p2}}\left(\frac{p-1}{p}\right)^{p\bigl(p+\tfrac12\bigr)}\E({X_n^*})^p\le \E S_n^p(X)\le \left(\frac2{p(p-1)}\right)^{\tfrac{p}2}\E (X_n^*)^p,
	 \end{equation}
	 or, equivalently,
	 \begin{equation}
	 \frac{1}{p\sqrt{2}}\left(\frac{p-1}{p}\right)^{p+\tfrac12}\bigl(\E({X_n^*})^p\bigr)^{\tfrac 1p}\le \bigl(\E S_n^p(X)\bigr)^{\tfrac 1p}\le\sqrt{\frac2{p(p-1)}}\bigl(\E (X_n^*)^p\bigr)^{\tfrac 1p}.
	 \end{equation}
	 	  The proof of \eqref{wz1} is complete.

\end{document}